\documentclass[paper=letter, fontsize=12pt]{scrartcl}

\usepackage{fourier,amssymb,epsfig,psfrag,amsmath,art,cite,wrapfig,multicol,color} 
\DeclareMathAlphabet{\mathcal}{OMS}{cmsy}{m}{n} 
\usepackage{subfigure}[1995/03/06] 
\usepackage[english]{babel} 

\newfont{\smalll}{cmr8}

\def\IR{\mathbb{R}}
\def\IS{\mathbb{S}}
\def\zero{\mathbf{0}}
\def\one{\mathbf{1}}

\def\diag{\mathrm{diag}}

\def\eye{\mathrm{I}}

\def\IS{\hbox{I\hskip-.1em S}}
\def\IC{\hbox{C\hskip-
.5em\raise.5ex\hbox{$\scriptscriptstyle\mid$}}\ }
\def\Ic{\hbox{\smalll C\hskip-
.5em\raise.3ex\hbox{$\scriptscriptstyle\mid$}}\ }

\def\T={\buildrel {\scriptscriptstyle\triangle} \over =}

\def\sqr#1#2{{\vcenter{\vbox{\hrule height.#2pt\hbox{\vrule
width.#2pt height#1pt \kern#1pt\vrule width.#2pt}\hrule
height.#2pt}}}}

\def\diag{\mathop{\rm diag}}

\def\block-diag{\mathop{\rm block{\scriptstyle -}diag}}

\def\pmbb#1{\setbox0=\hbox{#1}\raise 0.5ex\box0}

\def\norm#1{\|#1\|}

\newcommand{\bequ}{\begin{eqnarray}}
\newcommand{\eequ}{\end{eqnarray}}
\newcommand{\mT}{^\mathrm{T}}

\newcommand{\rom}{\mathrm}

\newcommand {\beq}      {\begin{equation}}
\newcommand {\eeq}      {\end{equation}}

\def\IR{{\mathbb R}}
\def\IC{{\mathbb C}}

\def\IS{{\mathbb S}}

\newcommand{\horrule}[1]{\rule{\linewidth}{#1}}

\begin{document}

\definecolor{tRed}{RGB}{250,25,0}
\newcommand{\red}[1]{\color{tRed}{#1}} 


\title{
\vspace{-0.3in} 	
\usefont{OT1}{bch}{b}{n}
\normalfont \normalsize \textsc{} \\ [25pt] 
\horrule{2pt} \\[0.4cm] 
\LARGE Control of Networked Multiagent Systems\\with Uncertain Graph Topologies$^\dagger${\thanks{\hspace{-0.15cm}$^\dagger$\hspace{0.08cm} This research was supported in part by the University of Missouri Research Board.\newline\indent \hspace{-0.844cm}$^*$\hspace{0.08cm} Corresponding author: 400 West, 13th Street, Rolla, MO 65409 (Address); +1 573 341 7702 (Phone); +1 573 341 6899 (Fax); yucelen@mst.edu (Email).}} \\
\horrule{2pt} \\[0.3cm]
}

\def\allAUTHORS
{\begin{tabular}{c}
   \textbf{Tansel~Yucelen}$^*$ and \textbf{John~Daniel~Peterson}\\
   \textit{Department of Mechanical and Aerospace Engineering}\\
   \textit{Missouri University of Science and Technology}\\
\\
\\
 \end{tabular}
\ \\[-1.5em]
\begin{tabular}{c}
    \textbf{Kevin L. Moore}\\
    \textit{Department of Electrical Engineering and Computer Science}\\
    \textit{Colorado School of Mines}\\
\\
\\
 \end{tabular}
}
\author{\allAUTHORS}

\date{January 2015} 

\maketitle 

\baselineskip 16pt


\vspace*{0cm} \vspace{-0.5em} 

\textbf{Abstract} --- Multiagent systems consist of agents that locally exchange information through a physical network subject to a graph topology. 
Current control methods for networked multiagent systems assume the knowledge of graph topologies in order to design distributed control laws for achieving desired global system behaviors. 
However, this assumption may not be valid for situations where graph topologies are subject to uncertainties either due to changes in the physical network or the presence of modeling errors especially for multiagent systems involving a large number of interacting agents. 
Motivating from this standpoint, this paper studies distributed control of networked multiagent systems with uncertain graph topologies. 
The proposed framework involves a controller architecture that has an ability to adapt its feedback gains in response to system variations. 
Specifically, we analytically show that the proposed controller drives the trajectories of a networked multiagent system subject to a graph topology with time-varying uncertainties to a close neighborhood of the trajectories of a given reference model having a desired graph topology. 
As a special case, we also show that a networked multiagent system subject to a graph topology with constant uncertainties asymptotically converges to the trajectories of a given reference model. 
Although the main result of this paper is presented in the context of average consensus problem, the proposed framework can be used for many other problems related to networked multiagent systems with uncertain graph topologies. 

\baselineskip 16pt


\vspace{0.75em}


\vspace{0.75em}

\clearpage \baselineskip=23.5pt \setcounter{page}{1}


\section{Introduction} 

Multiagent systems consist of agents that locally exchange information through a physical network subject to a graph topology. 
Current control methods for networked multiagent systems assume the knowledge of graph topologies in order to design distributed control laws for achieving desired global system behaviors (see, for example, \cite{mesbahi2010graph} and references therein). 
However, this assumption may not be valid for situations where graph topologies are subject to uncertainties either due to changes in the physical network or the presence of modeling errors especially for multiagent systems involving a large number of interacting agents. 

Uncertain nature of networked multiagent systems has received a considerable attention recently, including notable results \cite{sundaram2008distributed,sundaram2008distributed2,franceschelli2008motion,pasqualetti2007distributed,pasqualetti2012consensus,shames2011distributed,yucelen2012control,yucelen2013distributed,fazlyab2014robust}. 
For achieving desired multiagent system behavior, \cite{sundaram2008distributed,sundaram2008distributed2} make a specific assumption on the network connectivity other than the standard assumption on the connectedness of networked agents. 
The authors of \cite{franceschelli2008motion} excite the multiagent system in order to detect and isolate the uncertain agents from the network topology. 
Like \cite{sundaram2008distributed,sundaram2008distributed2}, a computationally expensive and not scalable algorithm is proposed in \cite{pasqualetti2007distributed,pasqualetti2012consensus} based on input observers technique, where the effect of uncertain agents on the overall multiagent system performance is quantified. 
An extension of this work is also given in \cite{shames2011distributed} that focuses on the detection and isolation of uncertain agents. 
The authors of \cite{yucelen2012control,yucelen2013distributed} use an adaptive control approach in order to suppress the effect of uncertain agents on the overall multiagent system performance without making specific assumptions on the fraction of misbehaving agents. 
A common similarity of the approaches documented in \cite{sundaram2008distributed,sundaram2008distributed2,franceschelli2008motion,pasqualetti2007distributed,pasqualetti2012consensus,shames2011distributed,yucelen2012control,yucelen2013distributed} is that they model uncertainties in the agent dynamics as additive perturbations that do not depend on the state of agents, where these results are not applicable to the networked multiagent systems with graph topology uncertainties since such uncertainties depend on the state of agents. 

One relevant work to the results of this paper is recently appeared in \cite{fazlyab2014robust}, where the authors utilize adaptive and sliding mode control methodologies in order to enforce a networked multiagent system subject to an uncertain graph topology to follow a given reference model having a desired graph topology. 
However, the result in \cite{fazlyab2014robust} may require a centralized information exchange among networked agents in general due to the structure of the proposed control algorithm (see (8) or (18) of \cite{fazlyab2014robust}). 
Other important results, which are related to this paper, are presented in \cite{zhou2007topology,tang2008adaptive,xu2009structure} without requiring a centralized information exchange. 
However, these results hold for graph topologies subject to constant uncertainties only. 

In this paper, we study distributed control of networked multiagent systems with uncertain graph topologies. 
The proposed framework involves a novel controller architecture that has an ability to adapt its feedback gains in response to system variations. 
Specifically, we analytically show that the proposed controller drives the trajectories of a networked multiagent system subject to a graph topology with time-varying uncertainties to a close neighborhood of the trajectories of a given reference model having a desired graph topology. 
As a special case, we also show that a networked multiagent system subject to a graph topology with constant uncertainties asymptotically converges to the trajectories of a given reference model. 
Although the main result of this paper is presented in the context of average consensus problem, the proposed framework can be used for many other problems related to networked multiagent systems with uncertain graph topologies. 


\section{Notation, Definitions, and Graph-Theoretic Notions} 

The notation used in this paper is fairly standard. 
Specifically, $\IR$ denotes the set of real numbers, 
$\IR^n$ denotes the set of $n \times 1$ real column vectors,
$\IR^{n \times m}$ denotes the set of $n \times m$ real matrices,
$\IR_+$ denotes the set of positive real numbers,
$\IR_+^{n \times n}$ (resp., $\overline{\IR}_+^{\hspace{0.1em} n \times n}$) denotes the set of $n \times n$ positive-definite (resp., nonnegative-definite) real matrices,
$\IS_+^{n \times n}$ (resp., $\overline{\IS}_+^{\hspace{0.1em} n \times n}$) denotes the set of $n \times n$ symmetric positive-definite (resp., symmetric nonnegative-definite) real matrices, 
$\zero_n$ (resp., $\one_n$ ) denotes the $n \times 1$ vector of all zeros (resp., ones), 
and $\eye_n$ denotes the $n \times n$ identity matrix. 
Furthermore, we write $(\cdot)\mT$ for transpose,
$\norm{\cdot}_2$ for the Euclidian norm, 
$\lambda_{\rom{min}}(A)$ (resp., $\lambda_{\rom{max}}(A)$) for the minimum (resp., maximum) eigenvalue of the Hermitian matrix $A$,
$\lambda_{i}(A)$ for the $i$-th eigenvalue of $A$ ($A$ is symmetric and the eigenvalues are ordered from least to greatest value),
$\diag(a)$ for the diagonal matrix with the vector $a$ on its diagonal, and 
$[A]_{ij}$ for the entry of the matrix $A$ on the $i$-th row and $j$-th column. 

Next, we recall some of the basic notions from graph theory, where we refer to \cite{mesbahi2010graph} for further details.
In the multiagent literature, graphs are broadly adopted to encode interactions in networked systems.
An \textit{undirected} graph $\mathcal{G}$ is defined by a set $\mathcal{V}_{\mathcal{G}}=\{1,\ldots,n\}$ of \textit{nodes}
and a set $\mathcal{E}_{\mathcal{G}} \subset \mathcal{V}_{\mathcal{G}} \times \mathcal{V}_{\mathcal{G}}$ of \textit{edges}.
If $(i,j) \in \mathcal{E}_{\mathcal{G}}$, then the nodes $i$ and $j$ are \textit{neighbors} and the neighboring relation is indicated with $i \sim j$.
The \textit{degree} of a node is given by the number of its neighbors.
Letting $d_i$ be the degree of node $i$, then the \textit{degree} matrix of a graph $\mathcal{G}$, $\mathcal{D}(\mathcal{G}) \in \IR^{n \times n}$, is given by $\mathcal{D}(\mathcal{G}) \triangleq \diag(d), \ d=[d_1,\ldots,d_n]\mT$.
A \textit{path} $i_0 i_1 \ldots i_L$ is a finite sequence of nodes such that $i_{k-1} \sim i_k$, $k=1, \ldots, L$,
and a graph $\mathcal{G}$ is \textit{connected} if there is a path between any pair of distinct nodes.
The \textit{adjacency} matrix of a graph $\mathcal{G}$, $\mathcal{A}(\mathcal{G}) \in \IR^{n \times n}$, is given by
\[
   [\mathcal{A}(\mathcal{G})]_{ij}
   \triangleq
   \left\{ \begin{array}{cl}
      1, &\mbox{ if $(i,j)\in\mathcal{E}_\mathcal{G}$},
      \\
      0, &\mbox{otherwise}.
   \end{array} \right.
   \label{AdjMat}
\]
The \textit{Laplacian} matrix of a graph, $\mathcal{L}(\mathcal{G}) \in \overline{\IS}_+^{\hspace{0.1em} n \times n}$, playing a central role in many graph theoretic treatments of multiagent systems, is given by $\mathcal{L}(\mathcal{G}) \triangleq \mathcal{D}(\mathcal{G}) - \mathcal{A}(\mathcal{G})$, where the spectrum of the Laplacian of a connected, undirected graph $\mathcal{G}$ can be ordered as
	\begin{equation}
		0 = \lambda_1(\mathcal{L}(\mathcal{G}))<\lambda_2(\mathcal{L}(\mathcal{G}))\le \cdots \le \lambda_n(\mathcal{L}(\mathcal{G})), \label{LapSpec}
	\end{equation}
with $\one_n$ as the eigenvector corresponding to the zero eigenvalue $\lambda_1(\mathcal{L}(\mathcal{G}))$ and $\mathcal{L}(\mathcal{G}) \one_n = \zero_n$ and $\rom{e}^{\mathcal{L}(\mathcal{G})} \one_n = \one_n$. 


\section{Problem Formulation} 

Consider a multiagent system consisting of $n$ agents that locally exchange information according to a connected, undirected \textit{uncertain} graph $\mathcal{G}_\mathrm{u}$ with nodes and edges representing agents and interagent information exchange links, respectively. 
We assume that the network is static, and hence, agent evolution will not cause edges to appear or disappear in the network. 
Specifically, let $x_i(t)\in\IR$ denote the state of node $i$ at time $t\in\overline{\IR}_+$, whose dynamics is described by
\bequ
\dot{x}_i(t)&=&-\alpha_i(t) x_i(t)+\sum_{i \sim j}\beta_{ij}(t) x_j(t)+u_i(t), \quad x_i(0)=x_{i0}, \label{eq:01}
\eequ
where $\alpha_i(t)\in\IR$ and $\beta_{ij}(t)\in\IR$ are unknown bounded coefficients of the graph $\mathcal{G}_\mathrm{u}$ with bounded time derivatives, and $u_i(t)\in\IR$, $t\in\overline{\IR}_+$, is the control input of node $i$. 
In this paper, we are interested to design a distributed control input $u_i(t)$, $t\in\overline{\IR}_+$, such that (\ref{eq:01}) achieves average consensus approximately (or asymptotically, i.e., $x(t) \rightarrow (\one_n \one_n\mT/n)x_0$ as $t \rightarrow \infty$, $x(t)=[x_1(t),\ldots,x_n(t)]\mT\in\IR^n$) in the presence of an uncertain graph topology.

\begin{figure}[t!] \hspace{2.5cm} \includegraphics[scale=0.45]{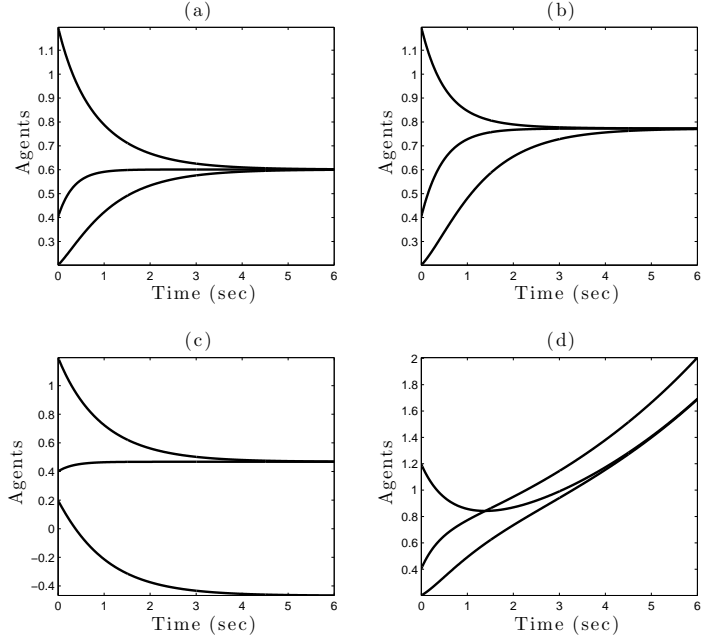} 
\caption{Trajectories of three agents on a line graph subject to initial conditions $(x_{10},x_{20},x_{30}=(0.2,0.4,1.2)$ for (a) $(\alpha_1,\alpha_2,\alpha_3)=(1,2,1)$ and $(\beta_{12},\beta_{21},\beta_{23},\beta_{32})=(1,1,1,1)$, (b) $(\alpha_1,\alpha_2,\alpha_3)=(1,1.1,1)$ and $(\beta_{12},\beta_{21},\beta_{23},\beta_{32})=(1,0.1,1,1)$, (c) $(\alpha_1,\alpha_2,\alpha_3)=(1,2,1)$ and $(\beta_{12},\beta_{21},\beta_{23},\beta_{32})=(-1,-1,1,1)$, and (d) $(\alpha_1,\alpha_2,\alpha_3)=(1,1.5,1)$ and $(\beta_{12},\beta_{21},\beta_{23},\beta_{32})=(1,1,1,1)$.}
\label{example:figure:01} 
\end{figure}

\textbf{Remark 1.} In the absence of proper control inputs $u_i(t)\in\IR$, $t\in\overline{\IR}_+$, (\ref{eq:01}) cannot necessarily achieve average consensus. 
To elucidate this point, let unknown coefficients of the graph $\mathcal{G}_\mathrm{u}$ be constant, i.e., $(\alpha_i(t),\beta_{ij}(t))=(\alpha_i,\beta_{ij})$, and consider four cases given in Figure \ref{example:figure:01} that show trajectories of three agents on a line graph subject to initial conditions $(x_{10},x_{20},x_{30})=(0.2,0.4,1.2)$. 
Since $\alpha_i=\sum_{i \sim j} \beta_{ij}$ and $\beta_{ij}=\beta_{ji}$ in case (a), this case results in average consensus at point $(\one_n \one_n\mT/n)x_0=0.6$. 
Since $\beta_{12}\neq\beta_{21}$ in case (b), this case does not result in average consensus at point $0.6$. 
Case (c) considers a multiagent system with antagonistic interactions \cite{altafini2013consensus}, and hence, it does not result in average consensus due to the existence of multiple equilibrium points. 
Finally, since $\alpha_i\neq\sum_{i \sim j} \beta_{ij}$ for $i=2$ in case (d), this case does not result in average consensus as well. 
In summary, (\ref{eq:01}) results in average consensus if $\alpha_i=\sum_{i \sim j} \beta_{ij}$ and $\beta_{ij}=\beta_{ji}\in\IR_+$ \cite{olfati2004consensus}. 
However, this cannot be justified due to unknown coefficients of the graph $\mathcal{G}_\mathrm{u}$, and hence, one needs to design proper control inputs $u_i(t)\in\IR$, $t\in\overline{\IR}_+$. 

Next, we propose a control input $u_i(t)\in\IR$, $t\in\overline{\IR}_+$, to drive the trajectories of (\ref{eq:01}) to a close neighborhood of a given reference model having a desired graph topology without requiring a centralized information exchange among networked agents. 
For this purpose, consider the reference model that locally exchange information according to a connected, undirected graph $\mathcal{G}$ given by
\bequ
\dot{r}_i(t)&=&-\sum_{i \sim j}\bigl(r_i(t)-r_j(t)\bigl), \quad r_i(0)=x_{i0}, \label{eq:02}
\eequ
where $r_i(t)\in\IR$, $t\in\overline{\IR}_+$, denotes the state of the reference model for node $i$. 
Note that 
\bequ
\lim_{t\rightarrow\infty}r(t)=(\one_n \one_n\mT/n)x_0, \label{extra:eq:01}
\eequ
where $r(t)=[r_1(t),\ldots,r_n(t)]\mT\in\IR^n$. 
Throughout this paper we assume that the nodes and edges of graphs $\mathcal{G}$ and $\mathcal{G}_\mathrm{u}$ coincide, however the graph $\mathcal{G}_\mathrm{u}$ is subject to unknown coefficients $\alpha_i(t)$ and $\beta_{ij}(t)$, as discussed earlier. 

\textbf{Remark 2.} The reference model given by (\ref{eq:02}) can be easily extended to
\bequ
\dot{r}_i(t)&=&-\sum_{i \sim j}\xi_{ij}\bigl(r_i(t)-r_j(t)\bigl), \quad r_i(0)=x_{i0}, \label{eq:03}
\eequ
without changing the following results of this paper, where $\xi_{ii}=\sum_{i \sim j}\xi_{ij}$ and $\xi_{ij}=\xi_{ji}\in\IR_+$. 

\textbf{Remark 3.} Note that (\ref{eq:02}) can be equivalently written as
\bequ
\dot{r}_i(t)&=&-d_i r_i(t)+\sum_{i \sim j}r_j(t), \quad r_i(0)=x_{i0}, \label{eq:04}
\eequ
where $d_i$ is the degree of node $i$ on graph $\mathcal{G}$. 
Therefore, if one knows the coefficients $\alpha_i(t)$ and $\beta_{ij}(t)$ of \ref{eq:01}, then the control input
\bequ
u_i(t)&=&-\bigl(d_i - \alpha_i(t)\bigl)x_i(t)+\sum_{i \sim j} \bigl(1-\beta_{ij}(t)\bigl)x_j(t), \label{eq:05}
\eequ
results in average consensus at point $(\one_n \one_n\mT/n)x_0$.

Since the control input (\ref{eq:05}) given in Remark 3 is not feasible due to unknown coefficients $\alpha_i(t)$ and $\beta_{ij}(t)$, we propose the adaptive control input given by
\bequ
 u_i(t)&=& -k_i \bigl(x_i(t)-r_i(t)\bigl)-\hat{w}_i(t)x_i(t)-\sum_{i \sim j} \hat{w}_{ij}(t)x_j(t), \label{eq:06}
\eequ
where $k_i\in\IR_+$ for at least one agent or a subset of agents (and $k_i=0$ for others), and the estimates $\hat{w}_i(t) \in\IR$ and $\hat{w}_{ij}(t)\in\IR$, $t\in\overline{\IR}_+$, satisfy the update laws 
\bequ
\dot{\hat{w}}_i(t) &=& \gamma_i \mathrm{Proj}\Bigl(\hat{w}_i(t), \ x_i(t) \bigl(x_i(t)-r_i(t)\bigl)\Bigl), \quad \hat{w}_i(0)=\hat{w}_{i0}, \label{eq:07} \\
\dot{\hat{w}}_{ij}(t) &=& \gamma_{ij} \mathrm{Proj}\Bigl(\hat{w}_{ij}(t), \ x_j(t) \bigl(x_i(t)-r_i(t)\bigl)\Bigl), \quad \hat{w}_{ij}(0)=\hat{w}_{ij0}, \label{eq:08}
\eequ
with $\gamma_i\in\IR_+$ and $\gamma_{ij}\in\IR_+$ being the corresponding learning rates. 
In the update laws given by (\ref{eq:07}) and (\ref{eq:08}), $\mathrm{Proj}$ denotes the projection operator \cite{pomet1992adaptive,yucelen2012robust}, which is used to keep the estimates $\hat{w}_i(t)$ and $\hat{w}_{ij}(t)$ bounded for all $t\in\overline{\IR}_+$. 
In the next section, we analytically show that the proposed adaptive control input given by (\ref{eq:06}) along with the update laws (\ref{eq:07}) and (\ref{eq:08}) drives the trajectories of (\ref{eq:01}) to a close neighborhood of the reference model trajectories given by (\ref{eq:03}). 


\section{Stability Analysis} 
In this section, we establish stability properties of the proposed adaptive control input given by (\ref{eq:06}) along with the update laws (\ref{eq:07}) and (\ref{eq:08}). 
For this purpose, let
\bequ
e_i(t) \triangleq x_i(t) - r_i(t), \quad t\in\overline{\IR}_+,
\eequ 
\vspace{0cm}
denote the local error dynamics that satisfy
\bequ
\hspace{-1.5cm}\dot{e}_i(t)&\hspace{-0.01cm}=\hspace{-0.01cm}& -\alpha_i(t) x_i(t)+\sum_{i \sim j}\beta_{ij}(t)x_j(t)+u_i(t)+d_i r_i(t)-\sum_{i \sim j}r_j(t) \nonumber\\
&\hspace{-0.01cm}=\hspace{-0.01cm}& -\alpha_i(t) x_i(t)+\sum_{i \sim j}\beta_{ij}(t)x_j(t)+u_i(t)+d_i r_i(t)-\sum_{i \sim j}r_j(t) + d_ix_i(t) \nonumber\\&& - d_ix_i(t)  +\sum_{i \sim j} x_j(t) - \sum_{i\sim j}x_j(t)\nonumber\\
&\hspace{-0.01cm}=\hspace{-0.01cm}& -d_i e_i(t) +\sum_{i \sim j}e_j(t) + \bigl(d_i-\alpha_i(t)\bigl)x_i(t)+\sum_{i \sim j}\bigl(\beta_{ij}(t)-1\bigl)x_j(t)+u_i(t) \nonumber\eequ\bequ
\hspace{-0.2cm}&\hspace{-0.01cm}=\hspace{-0.01cm}& -\sum_{i \sim j}\bigl(e_i(t)-e_j(t)\bigl) + w_i(t) x_i(t)+\sum_{i \sim j}w_{ij}(t) x_j(t)+u_i(t) \nonumber\\
&\hspace{-0.01cm}=\hspace{-0.01cm}& -\sum_{i \sim j}\bigl(e_i(t)-e_j(t)\bigl) + w_i(t) x_i(t)+\sum_{i \sim j}w_{ij}(t) x_j(t)-k_i e_i(t)\nonumber\\&&
-\hat{w}_i(t)x_i(t)-\sum_{i \sim j} \hat{w}_{ij}(t)x_j(t) \nonumber\\
&\hspace{-0.01cm}=\hspace{-0.01cm}& -k_i e_i(t)-\sum_{i \sim j}\bigl(e_i(t)-e_j(t)\bigl) - \tilde{w}_i(t) x_i(t)-\sum_{i \sim j}\tilde{w}_{ij}(t) x_j(t), \quad e_i(0)=0, \label{stab:01}
\eequ
where 
\bequ
\tilde{w}_i(t)&\triangleq& \hat{w}_i(t)-w_i(t), \quad t\in\overline{\IR}_+,\\
\tilde{w}_{ij}(t)&\triangleq& \hat{w}_{ij}(t)-w_{ij}(t), \quad t\in\overline{\IR}_+,
\eequ
$w_i(t)\triangleq d_i-\alpha_i(t)$, $t\in\overline{\IR}_+$, and $w_{ij}(t)\triangleq\beta_{ij}(t)-1$, $t\in\overline{\IR}_+$.
In addition, it follows from (\ref{eq:07}) and (\ref{eq:08}) that
\bequ
\dot{\tilde{w}}_i(t) &=& \gamma_i \mathrm{Proj}\Bigl(\hat{w}_i(t), \ x_i(t) \bigl(x_i(t)-r_i(t)\bigl)\Bigl)-\dot{w}_i(t), \quad \tilde{w}_i(0)=\tilde{w}_{i0}, \label{stab:02} \\
\dot{\tilde{w}}_{ij}(t) &=& \gamma_{ij} \mathrm{Proj}\Bigl(\hat{w}_{ij}(t), \ x_j(t) \bigl(x_i(t)-r_i(t)\bigl)\Bigl)-\dot{w}_{ij}(t), \quad \tilde{w}_{ij}(0)=\tilde{w}_{ij0}, \label{stab:03}
\eequ
where $\tilde{w}_{i0}\triangleq\hat{w}_{i0}-w_i$ and $\tilde{w}_{ij0}\triangleq\hat{w}_{ij0}-w_{ij}$. 
Note that $\dot{w}_i(t)$ and $\dot{w}_{ij}(t)$ are bounded since it is assumed that unknown bounded coefficients $\alpha_i(t)$ and $\beta_{ij}(t)$ have bounded time derivates. 
We now state the following lemma necessary for the results of this section.

\textbf{Lemma 1.}
	Let $K=\rom{diag}(k)$, $k=[k_1, k_2, \ldots, k_n]\mT$, $k_i \in \overline{\mathbb{R}}_+$, $i=1,\ldots,n$, and assume that at least one element of $k$ is nonzero. 
	Then, for the Laplacian of a connected, undirected graph, 
	\bequ
	\mathcal{F}(\mathcal{G})\triangleq\mathcal{L}(\mathcal{G})+K\in\IS_+^{n \times n}, \label{matrix_f}
	\eequ
	and $\rom{det}(\mathcal{F}(\mathcal{G}))\neq0$.
	
\textit{Proof.} 
Consider the decomposition $K=K_1+K_2$, where $K_1\triangleq\rom{diag}([0, \ldots, 0, \phi_i, 0, \ldots, 0]\mT)$ and $K_2\triangleq K-K_1$, where $\phi_i$ denotes the smallest nonzero diagonal element of $K$ appearing on its $i$-th diagonal, so that $K_2\in\overline{\IS}_+^{n \times n}$. 
From the Rayleigh's Quotient \cite{Lay2006}, the minimum eigenvalue of $\mathcal{L}(\mathcal{G})+K_1$ can be given by
\bequ
\lambda_{\text{min}}(\mathcal{L}(\mathcal{G})+K_1) &\hspace{-0.1cm}=\hspace{-0.1cm}& \min\limits_{x} \{x\mT\bigl(\mathcal{L}(\mathcal{G})+K_1\bigl)x \: | \: x\mT x = 1\}, \label{added:equation:1}
\eequ
where $x$ is the eigenvector corresponding to this minimum eigenvalue. 
Note that since $\mathcal{L}(\mathcal{G})\in\overline{\IS}_+^{n \times n}$ and $K_1\in\overline{\IS}_+^{n \times n}$, and hence, $\mathcal{L}(\mathcal{G})+K_1$ is real and symmetric, $x$ is a real eigenvector. 
Now, expanding (\ref{added:equation:1}) as
\bequ
		 x\mT\left( \mathcal{L}(\mathcal{G}) + K_1\right)x &=& \sum\limits_{i \thicksim j}(x_i-x_j)^2 + \phi_i x_i^2,
		\label{eq:lem2}
\eequ
and noting that the right hand side of (\ref{eq:lem2}) is zero only if $x\equiv0$, it follows that $\lambda_{\text{min}}(\mathcal{L}(\mathcal{G}) + K_1) > 0$, and hence, $\mathcal{L}(\mathcal{G}) + K_1\in\IS_+^{n \times n}$. 
Finally, let $\lambda$ be an eigenvalue of $\mathcal{F}(\mathcal{G})=\mathcal{L}(\mathcal{G}) + K_1 + K_2$. 
Since $\lambda_\rom{min}(\mathcal{L}(\mathcal{G}) + K_1)>0$ and $\lambda_\rom{min}(K_2)=0$, it follows from Fact 5.11.3 of \cite{bernstein2009matrix} that $\lambda_\rom{min}(\mathcal{L}(\mathcal{G}) + K_1) + \lambda_\rom{min}(K_2) \le \lambda$, 
and hence, $\lambda>0$, which implies that (\ref{matrix_f}) holds and $\rom{det}(\mathcal{F}(\mathcal{G}))\neq0$. \hfill $\blacksquare$

The next theorem presents the first result of this section.

\textbf{Theorem 1.} 
Consider the networked multiagent system given by (\ref{eq:01}) subject to an uncertain graph topology, the reference model given by (\ref{eq:02}), the adaptive control input given by (\ref{eq:06}), and the update laws given by (\ref{eq:07}) and (\ref{eq:08}). 
Then, the solution $\bigl(e_i(t),\tilde{w}_i(t),\tilde{w}_{ij}(t)\bigl)$ of the closed-loop dynamical system given by (\ref{stab:01}), (\ref{stab:02}), and (\ref{stab:03}) is bounded for all $\bigl(0,\tilde{w}_{i0},\tilde{w}_{ij0}\bigl)$, $t\in\overline{\IR}_+$, and $(i,j)$. 

\textit{Proof.} 
First consider the quadratic function given by
\bequ
V_i(e_i,\tilde{w}_i,\tilde{w}_{ij})&=&\frac{1}{2}\Bigl( e_i^2 + \gamma_i^{-1} \tilde{w}_i^2 + \sum_{i \sim j}\gamma_{ij}^{-1}\tilde{w}_{ij}^2 \Bigl), \label{proof:01}
\eequ
and note that $V_i(0,0,0)=0$ and $V_i(e_i,\tilde{w}_i,\tilde{w}_{ij})\in\IR_+$, $(e_i,\tilde{w}_i,\tilde{w}_{ij})\neq(0,0,0)$. 
Furthermore, $V_i(e_i,$ $\tilde{w}_i,\tilde{w}_{ij})$ is radially unbounded. 
Differentiating (\ref{proof:01}) along the closed-loop trajectories of (\ref{stab:01}), (\ref{stab:02}), and (\ref{stab:03}) yields
\bequ
\dot{V}_i\bigl(e_i(t),\tilde{w}_i(t),\tilde{w}_{ij}(t)\bigl)&\le&-e_i(t)\sum_{i \sim j}\bigl(e_i(t)-e_j(t)\bigl)-k_i e_i^2(t)+w_i^*, \label{proof:02}
\eequ
where $w_i^*$ is an upper bound satisfying $\bigl{|}\bigl{|}\gamma_i^{-1}\bigl(\hat{w}_i(t)-w_i(t)\bigl)\dot{w}_i(t)+\sum_{i \sim j}\gamma_{ij}^{-1}\bigl(\hat{w}_{ij}(t)-w_{ij}(t)\bigl)$    $\cdot\dot{w}_{ij}(t)  \bigl{|}\bigl{|}_2 \le w_i^*$, $t\in\overline{\IR}_+$.
Note that $w_i^*$ exists since all the terms inside the norm operator are bounded and projection operator is used for the estimates $\hat{w}_i(t)$ and $\hat{w}_{ij}(t)$. 
Now, consider the Lyapunov function candidate given by
\bequ
V(\cdot)&=&\sum_{i=1}^n V_i(e_i,\tilde{w}_i,\tilde{w}_{ij}). \label{proof:03}
\eequ
The time derivative of (\ref{proof:03}) is given using (\ref{proof:02}) as 
\bequ
\dot{V}(\cdot) &\le& -e\mT(t) \bigl( \mathcal{L}(\mathcal{G}) + K\bigl)e(t) + w^*, \quad w^*\triangleq \sum_{i=1}^n w_i^*, \label{proof:04}
\eequ
where $\mathcal{L}(\mathcal{G})$ denotes the Laplacian matrix of (\ref{eq:02}), $K\triangleq\rom{diag}(k)$, $k=[k_1, k_2, \ldots, k_n]\mT$, $k_i \in \overline{\mathbb{R}}_+$, and $e(t)=[e_1(t),\ldots,e_n(t)]\mT$. 
From the definition of the adaptive control input in (\ref{eq:06}), notice that at least one element of $k$ is nonzero. 
This implies from Lemma 1 that $\mathcal{L}(\mathcal{G}) + K\in\IS_+^{n \times n}$ and $\rom{det}(\mathcal{L}(\mathcal{G}) + K)\neq0$, and hence, $\lambda_\rom{min}\bigl( \mathcal{L}(\mathcal{G}) + K\bigl) \norm{e(t)}_2 \le e\mT(t) \bigl( \mathcal{L}(\mathcal{G}) + K\bigl)e(t)$. 
Now, since $\dot{V}(\cdot)\le0$ when $\norm{e(t)}_2\ge w^*/\lambda_\rom{min}\bigl( \mathcal{L}(\mathcal{G}) + K\bigl)$, it follows that the closed-loop dynamical system given by (\ref{stab:01}), (\ref{stab:02}), and (\ref{stab:03}) is bounded for all $\bigl(0,\tilde{w}_{i0},\tilde{w}_{ij0}\bigl)$, $t\in\overline{\IR}_+$, and $(i,j)$. \hfill $\blacksquare$

\textbf{Remark 4.} In order to drive the trajectories of (\ref{eq:01}) to a close neighborhood of the reference model trajectories given by (\ref{eq:03}), the perturbation term $w^*$ in (\ref{proof:04}) needs to be small. 
This holds if the time derivative of unknown coefficients $\alpha_i(t)$ and $\beta_{ij}(t)$ is small. 
If this is not true, then one can increase the learning rates $\gamma_i$ and $\gamma_{ij}$ to make $w^*$ small. 

As a special case when the unknown coefficients are constant, i.e., $(\alpha_i(t),\beta_{ij}(t))=(\alpha_i,\beta_{ij})$, the next theorem shows that the proposed adaptive control input given by (\ref{eq:06}) along with the update laws (\ref{eq:07}) and (\ref{eq:08}) asymptotically drives the trajectories of (\ref{eq:01}) to the reference model trajectories given by (\ref{eq:03}). 

\textbf{Theorem 2.} 
Consider the networked multiagent system given by (\ref{eq:01}) subject to an uncertain graph topology, the reference model given by (\ref{eq:02}), the adaptive control input given by (\ref{eq:06}), and the update laws given by (\ref{eq:07}) and (\ref{eq:08}). 
Then, the solution $\bigl(e_i(t),\tilde{w}_i(t),\tilde{w}_{ij}(t)\bigl)$ of the closed-loop dynamical system given by (\ref{stab:01}), (\ref{stab:02}), and (\ref{stab:03}) is Lyapunov stable for all $\bigl(0,\tilde{w}_{i0},\tilde{w}_{ij0}\bigl)$, $t\in\overline{\IR}_+$, and $(i,j)$, and $\lim_{t\rightarrow\infty}e_i(t)=0$ for all $i$. 
In addition, $\lim_{t\rightarrow\infty}x(t)=(\one_n \one_n\mT/n)x_0$. 

\textit{Proof.} 
To show the Lyapunov stability of the closed-loop dynamical system given by (\ref{stab:01}), (\ref{stab:02}), and (\ref{stab:03}), first consider the quadratic function given by (\ref{proof:01}). 
Differentiating (\ref{proof:01}) along the closed-loop trajectories of (\ref{stab:01}), (\ref{stab:02}), and (\ref{stab:03}) yields
\bequ
\dot{V}_i\bigl(e_i(t),\tilde{w}_i(t),\tilde{w}_{ij}(t)\bigl)&\le&-e_i(t)\sum_{i \sim j}\bigl(e_i(t)-e_j(t)\bigl)-k_i e_i^2(t). \label{proof:02b}
\eequ
Now, consider the Lyapunov function candidate given by (\ref{proof:03}), where the time derivative of (\ref{proof:03}) is given using (\ref{proof:02b}) as 
\bequ
\dot{V}(\cdot) &\le& -e\mT(t) \bigl( \mathcal{L}(\mathcal{G}) + K\bigl)e(t), \label{proof:04}
\eequ
where $\mathcal{L}(\mathcal{G})$ denotes the Laplacian matrix of (\ref{eq:02}), $K\triangleq\rom{diag}(k)$, $k=[k_1, k_2, \ldots, k_n]\mT$, $k_i \in \overline{\mathbb{R}}_+$, and $e(t)=[e_1(t),\ldots,e_n(t)]\mT$. 
From the definition of the adaptive control input in (\ref{eq:06}), notice that at least one element of $k$ is nonzero. 
This implies from Lemma 1 that $\mathcal{L}(\mathcal{G}) + K\in\IS_+^{n \times n}$ and $\rom{det}(\mathcal{L}(\mathcal{G}) + K)\neq0$. 
Hence, the closed-loop dynamical system given by (\ref{stab:01}), (\ref{stab:02}), and (\ref{stab:03}) is Lyapunov stable for all $\bigl(0,\tilde{w}_{i0},\tilde{w}_{ij0}\bigl)$, $t\in\overline{\IR}_+$, and $(i,j)$. 

Next, it follows from \cite{haddad2008nonlinear} that $\lim_{t\rightarrow\infty} e\mT(t) \bigl( \mathcal{L}(\mathcal{G}) + K\bigl)e(t)=0$ holds, which implies that $\lim_{t\rightarrow\infty}$ $e(t)=0$ as a consequence of $\rom{det}(\mathcal{L}(\mathcal{G}) + K)\neq0$. 
The result $\lim_{t\rightarrow\infty}x(t)=(\one_n \one_n\mT/n)x_0$ is now immediate. \hfill $\blacksquare$

\begin{figure}[t!] \hspace{2.5cm} \includegraphics[scale=0.45]{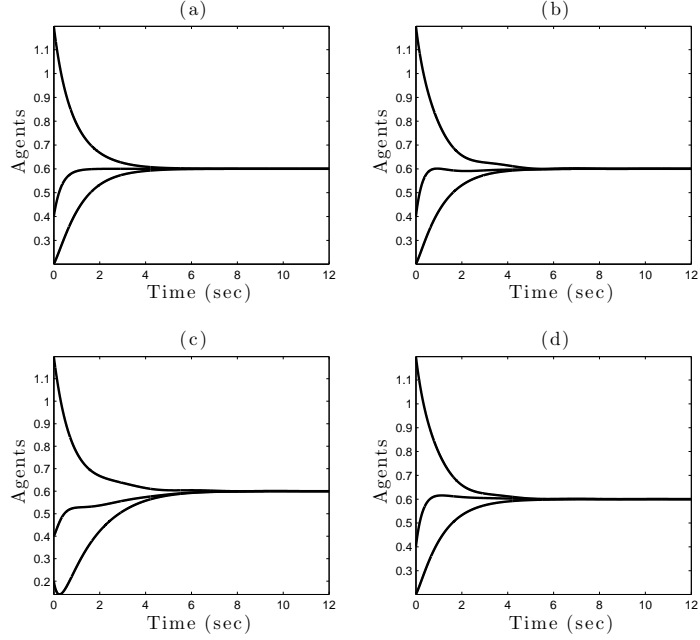} 
\caption{Trajectories of three agents on a line graph subject to initial conditions $(x_{10},x_{20},x_{30}=(0.2,0.4,1.2)$ and the proposed approach for (a) $(\alpha_1,\alpha_2,\alpha_3)=(1,2,1)$ and $(\beta_{12},\beta_{21},\beta_{23},\beta_{32})=(1,1,1,1)$, (b) $(\alpha_1,\alpha_2,\alpha_3)=(1,1.1,1)$ and $(\beta_{12},\beta_{21},\beta_{23},\beta_{32})=(1,0.1,1,1)$, (c) $(\alpha_1,\alpha_2,\alpha_3)=(1,2,1)$ and $(\beta_{12},\beta_{21},\beta_{23},\beta_{32})=(-1,-1,1,1)$, and (d) $(\alpha_1,\alpha_2,\alpha_3)=(1,1.5,1)$ and $(\beta_{12},\beta_{21},\beta_{23},\beta_{32})=(1,1,1,1)$.}
\label{example:figure:02} 
\end{figure}

\textbf{Remark 5.} 
We now revisit the example in Remark 1 and use the adaptive control input given by (\ref{eq:06}) with the update laws given by (\ref{eq:07}) and (\ref{eq:08}). 
In particular, we set $(k_1,k_2,k_3)=(5,5,0)$ and $\gamma_i=\gamma_{ij}=5$ for all $(i,j)$, and use zero initial conditions for the update laws. 
Figure \ref{example:figure:02} shows that the proposed approach achieves average consensus for all cases, i.e., $\lim_{t \rightarrow \infty}x_i(t)=0.6$ for all $i$, as expected from Theorem 2.

\section{Conclusion} 

In order to contribute to the previous studies in networked multiagent systems, we investigated an adaptive control methodology that has an ability to drive the trajectories of an uncertain and time-varying multiagent system to a close neighborhood of the trajectories of a given reference model having a desired graph topology. 
In the context of average consensus problem, we rigorously analyzed stability properties of this methodology using results from nonlinear systems theory and matrix mathematics. 
In addition, as a special case when the unknown coefficients of the graph are constant, we showed that the uncertain multiagent system asymptotically converges to the given reference model with the proposed control methodology. 
Illustrative examples indicated that the presented theory and its numerical results are compatible.  


\bibliographystyle{IEEEtran} \baselineskip 16pt
\bibliography{references.bib}
\end{document}